\documentclass[a4paper,11pt]{my-ijamas}

\usepackage{graphicx}
\usepackage{url}

\begin{document}


\title{Numerical approximation of the thermistor problem\thanks{Research Report CM07/I-16. The original publication is available at
http://www.isder.ceser.res.in/ijms.html,
Int. J. Math. Stat.,
Vol.~2, Nr.~S08, 2008, pp.~106--114.}}

\author{\textbf{Moulay Rchid Sidi Ammi$^1$ and Delfim F. M. Torres$^2$}}

\date{$^1$Department of Mathematics\\
University of Aveiro\\
3810-193 Aveiro, Portugal\\
\url{sidiammi@mat.ua.pt}\\ [0.3cm]
$^2$Department of Mathematics\\
University of Aveiro\\
3810-193 Aveiro, Portugal\\
\url{delfim@ua.pt}}

\maketitle


\begin{abstract}
\noindent \emph{We use a finite element approach based on Galerkin
method to obtain approximate steady state solutions of the
thermistor problem with temperature dependent electrical
conductivity.}

\medskip


\medskip


\end{abstract}


\section{Introduction}

In this paper we develop a method to approximate steady-state
solutions of the following one-dimensional thermistor problem:
\begin{equation} \label{eq1}
\frac{\partial u}{\partial t} - \frac{\partial}{\partial x}\left(
k(u)u_{x} \right )= \sigma (u)|\varphi_{x}|^{2}, \, \quad 0< x < 1,
\, \quad  t>0,
\end{equation}
subject to boundary and initial conditions,
\begin{equation} \label{eq2}
k(u)\frac{\partial u}{\partial x} = -\beta u \, \quad \mbox { on }
\partial \Omega \times (0, T),
\end{equation}
\begin{equation} \label{eq3}
u(x, 0)= 0, \, \quad 0\leq x \leq 1,
\end{equation}
and coupled with the electric potential equation:
\begin{equation} \label{eq4}
\left( \sigma(u) \varphi_{x} \right)_{x}=0,  \, \quad 0< x < 1, \,
\quad t>0,
\end{equation}
\begin{equation} \label{eq5}
 \frac{\partial \varphi}{\partial x}=
\overline{\varphi}(x, t) \, \quad \mbox{ on } \partial \Omega \, ,
\end{equation}
\begin{equation} \label{eq6}
\varphi(x, 0)= x, \, \quad 0\leq x \leq 1 \, .
\end{equation}

The motivation for studying this kind of problem is that
\eqref{eq1}-\eqref{eq6} has important implications for a variety
of technological processes. For example, it arises in the
analytical study of phenomena associated with the occurrence of
shear band in metal being deformed at high strain rates \cite{bl};
in the theory of gravitational equilibrium of polytropic stars
\cite{kn}; in the investigation of the fully turbulent behavior of
flows \cite{clmp}; in modelling aggregation of cells via
interaction with a chemical substance (chemotaxis) \cite{w}; and
specially in modelling electrical heating in a conductor
\cite{xu}. In this case, $u$ is the temperature of the conductor,
$\varphi$ is the electrical potential. Functions $\sigma(u)$ and
$k(u)$ are, respectively, the electrical and thermal
conductivities; $\beta$ is the heat transfer coefficient. The
condition \eqref{eq3} is a condition of Robin-type. When $\beta
=0$ it is called an adiabatic condition. Equation \eqref{eq1}
consists in the heat equation with Joule heating as a source;
\eqref{eq4} describes conservation of current in the conductor.

The thermistor problem has been extensively studied by several
authors \cite{ac,ci,se1,se2}, where existence and uniqueness of
solutions are given. Theoretical analysis, consisting in existence
of solutions with the required regularity and which ensure error
estimates of optimal order of convergence, are done in
\cite{cesl}. To construct a numerical approximation of the steady
state solution we use a numerical method to approximate the
solution of the parabolic problem. This approach has been used by
\cite{ba,ku} in the one-dimensional thermistor problem. Further,
in these last works authors consider the thermal conductivity $k$
equal to $1$ and a particular electrical conductivity, then they
obtain the exact solution $(\varphi(x,t)=x)$ of the conservation
problem \eqref{eq4}-\eqref{eq6} and so system
\eqref{eq1}-\eqref{eq6} of thermistor problem is reduced to the
following single heat conduction problem:
\begin{equation*}
\frac{\partial u}{\partial  t}- \frac{\partial^{2}u}{\partial^{2}
x}= \sigma (u),
\end{equation*}
subject to the boundary conditions \eqref{eq2}-\eqref{eq3}. In
this paper, we propose to solve both equations \eqref{eq1} and
\eqref{eq6} at the same time by using a finite element method and
a fully Crank-Nicolson approach. The formulation of the finite
element method is standard and is based on a variational
formulation of the continuous problem. In Section~\ref{sec2} we
give the variational formulation of problem
\eqref{eq1}-\eqref{eq6}. An algorithm for solving the problem is
then proposed in Section~\ref{sec3}. In Section~\ref{sec4},
numerical results are obtained for an appropriate test-problem.


\section{Variational formulation of the problem}
\label{sec2}

We divide the interval $\Omega = [0, 1]$ into $ N$ equal finite
elements $0=x_{0} < x_{1} < \ldots < x_{N}=1$. Let $(x_{j},
x_{j+1})$ be a partition of $\Omega$ and  $x_{j+1}-x_{j}=h=
\frac{1}{N}$  the step length. By $S$ we denote a basis of  the
usual pyramid functions:
\begin{equation*}
v_{j}=
\begin{cases}
\frac{1}{h}x+(1-j)
 & \mbox{ on } [x_{j-1}, x_{j}],\\
- \frac{1}{h}x+(1+j)
 & \mbox{ on } [x_{j}, x_{j+1}],\\
 0 & otherwise.
\end{cases}
\end{equation*}
As indicated above, it is convenient to proceed in two steps with
the derivation and analysis of the approximate solution of
\eqref{eq1}-\eqref{eq6}. First, we write the problem in weak or
variational form. We multiply the parabolic equation by $v_{j}$
(for $j$ fixed), integrate over $(0, 1)$, and apply Green's
formula on the left-hand side, to obtain
\begin{equation*}
\int_{\Omega} \frac{\partial u}{\partial t} v_{j} \, dx +
\int_{\Omega} k(u) \nabla u \nabla v_{j} \, dx - \int_{\partial
\Omega} k(u)\frac{\partial u}{\partial \nu} v_{j} \, ds =
\int_{\Omega} \sigma (u) |\nabla \varphi|^{2} v_{j} \, dx.
\end{equation*}
Using the boundary condition we get
\begin{equation}
\label{eq7}
\int_{\Omega} \frac{\partial u}{\partial t} v_{j} \, dx +
\int_{\Omega} k(u) \nabla u \nabla v_{j} \, dx + \int_{\partial
\Omega}\beta  u v_{j} \, ds = \int_{\Omega} \sigma (u) |\nabla
\varphi|^{2} v_{j} \, dx.
\end{equation}
We also have
\begin{equation}
\label{eq8}
\int_{\Omega}  \sigma(u) \nabla \varphi  \nabla v_{j} \, dx =
\int_{\partial \Omega} \sigma (u)\frac{\partial \varphi}{\partial
\nu} v_{j} \, ds.
\end{equation}
We now turn our attention to solve this system by discretization
with respect to the time variable. We introduce a time step $\tau$
and time levels $t_{n}= n\tau$, where $n$ is a nonnegative
integer, and denote by $u^{n}$ the approximation of $u(t_{n})$ to
be determined. We use the backward Euler Galerkin method, which is
defined by replacing the time derivative in \eqref{eq7} by a
backward difference $\frac{u^{n+1}-u^{n}}{\tau}$. So the
approximations $u^{n+1}$, $\varphi^{n+1}$ admit unique
representations
$$
u^{n+1}= \sum_{i=-1}^{N}\alpha_{i}^{n+1} v_{i} \, \quad
\varphi^{n+1}= \sum_{i=-1}^{N}\mu_{i}^{n+1} v_{i},
$$
where $\alpha_{i}^{n+1}$, $\mu_{i}^{n+1}$ are unknown real
coefficients to be determined. Then, after decoupling, we have
that
\begin{equation}
\label{eq9}
\int_{\Omega} \frac{u^{n+1}-u^{n}}{\tau} v_{j} \, dx + \int_{\Omega}
k(u^{n}) \nabla u^{n+1} \nabla v_{j} \, dx +\int_{\partial \Omega}
\beta u^{n+1} v_{j} \, ds =
  \int_{\Omega} \sigma (u^{n})
|\nabla \varphi^{n}|^{2} v_{j} \, dx,
\end{equation}
and
\begin{equation}\label{eq10}
\int_{\Omega}  \sigma(u^{n}) \nabla \varphi^{n+1}  \nabla v_{j} \,
dx = \int_{\partial \Omega} \sigma (u^{n})\frac{\partial
\varphi^{n+1}}{\partial \nu} v_{j} \, ds.
\end{equation}


\section{Formulation of the numerical method}
\label{sec3}
For scheme \eqref{eq10}, we have
\begin{equation*}
\begin{split}
& \sum_{i=1}^{N} \int_{\Omega} \sigma(u^{n}) \frac{\partial
v_{i}}{\partial x}  \frac{\partial v_{j}}{\partial x} \, dx\\
& = \mu_{j-1}^{n+1}\int_{x_{j-1}}^{x_{j}}\sigma(u^{n})
\frac{\partial v_{j-1}}{\partial x}  \frac{\partial v_{j}}{\partial
x} \, dx
 +\mu_{j}^{n+1}\int_{x_{j-1}}^{x_{j+1}}\sigma(u^{n})
(\frac{\partial v_{j}}{\partial x})^{2}   \, dx
 + \mu_{j+1}^{n+1}\int_{x_{j}}^{x_{j+1}}\sigma(u^{n})
\frac{\partial v_{j}}{\partial x}  \frac{\partial v_{j+1}}{\partial
x} \, dx \\
&=-\frac{\mu_{j-1}^{n+1}}{h^{2}} \int_{x_{j-1}}^{x_{j}}\sigma(u^{n})
\, dx +\frac{\mu_{j}^{n+1}}{h^{2}}
\int_{x_{j-1}}^{x_{j+1}}\sigma(u^{n}) \, dx
 - \frac{\mu_{j+1}^{n+1}}{h^{2}} \int_{x_{j}}^{x_{j+1}}\sigma(u^{n})
\, dx\\
& \simeq -\frac{\mu_{j-1}^{n+1}}{2h}\left( \sigma(u^{n}(x_{j}))+
\sigma(u^{n}(x_{j-1}))\right)
 + \frac{\mu_{j}^{n+1}}{h}\left( \sigma(u^{n}(x_{j+1}))+
\sigma(u^{n}(x_{j-1}))\right)\\
 & \qquad -\frac{\mu_{j+1}^{n+1}}{2h}\left( \sigma(u^{n}(x_{j+1}))+
\sigma(u^{n}(x_{j}))\right)\\
& \simeq -\frac{\mu_{j-1}^{n+1}}{2h}\left( \sigma(\alpha_{j}^{n})+
\sigma(\alpha_{j-1}^{n})\right)
 + \frac{\mu_{j}^{n+1}}{h}\left( \sigma(\alpha_{j+1}^{n})+
\sigma(\alpha_{j-1}^{n})\right)
 -\frac{\mu_{j+1}^{n+1}}{2h}\left( \sigma(\alpha_{j+1}^{n})+
\sigma(\alpha_{j}^{n})\right).
\end{split}
\end{equation*}
On the other hand, we have
\begin{equation*}
\begin{split}
 \int_{\partial \Omega}\sigma(u^{n})\frac{\partial
\varphi^{n+1}}{\partial \nu} v_{j} \, ds &= \int_{\partial
\Omega}\sigma(u^{n}) \overline{\varphi} v_{j} \, ds\\
& = \sigma(u^{n}(1)) \overline{\varphi}(1) v_{j}(1) -
\sigma(u^{n}(0)) \overline{\varphi}(0) v_{j}(0)\\
&=\begin{cases}
-\sigma(\alpha_{0}^{n})\overline{\varphi}(0) & \mbox{ if } j= 0,\\
 0 & \mbox{ if } j= 1, \ldots N-2, \\
 0 & \mbox{ if } j= N-1.
 \end{cases}
\end{split}
\end{equation*}
Using boundary conditions \eqref{eq2} and initial condition
\eqref{eq3}, it follows that
\begin{equation*}
\begin{split}
& \mu_{-1}^{n+1} = \mu_{1}^{n+1}- \mu_{0}^{n+1}-h
\overline{\varphi}(0),\\
 & \mu_{N}^{n+1} = h \overline{\varphi}(1)+ \mu_{N-1}^{n+1},\\
 &\alpha_{-1}^{n}= \alpha_{1}^{n}+ \left(
\frac{h\beta}{k(\alpha_{0}^{n-1})}-1 \right) \alpha_{0}^{n}.
 \end{split}
\end{equation*}
Then, we have the resulting system of equations:

for $j=0$,
\begin{multline}
\label{eq11}
 \left ( \sigma(\alpha_{0}^{n})+3 \sigma(\alpha_{-1}^{n}) + 2
\sigma(\alpha_{1}^{n}) \right) \mu_{0}^{n+1} -
  \left ( \sigma(\alpha_{-1}^{n})+2 \sigma(\alpha_{0}^{n}) +
\sigma(\alpha_{1}^{n}) \right) \mu_{1}^{n+1}\\
= - h \overline{\varphi}(0)( 3\sigma(\alpha_{0}^{n}) +
\sigma(\alpha_{-1}^{n})) \, ;
\end{multline}
for $j=1, \ldots , N-2$,
\begin{multline}
\label{eq12}
\begin{split}
&-\mu_{j-1}^{n+1}\left(\sigma(\alpha_{j}^{n})+
\sigma(\alpha_{j-1}^{n})\right)
 + 2\mu_{j}^{n+1}\left(\sigma(\alpha_{j+1}^{n})+\sigma(\alpha_{j-1}^{n})\right)
 + \mu_{j+1}^{n+1}\left(\sigma(\alpha_{j+1}^{n})+
\sigma(\alpha_{j}^{n})\right)= 0 \, ;
\end{split}
\end{multline}
for $j= N-1$,
\begin{multline}
\label{eq13}
 -\left (\sigma(\alpha_{N-1}^{n})+\sigma(\alpha_{N-2}^{n})\right )
\mu_{N-2}^{n+1} +  \left
(2\sigma(\alpha_{N-2}^{n})+\sigma(\alpha_{N}^{n})-
\sigma(\alpha_{N-1}^{n}) \right ) \mu_{N-1}^{n+1}\\
= h \overline{\varphi}(1)\left (\sigma(\alpha_{N}^{n})+
\sigma(\alpha_{N-1}^{n}) \right ).
\end{multline}
Coming back to \eqref{eq9}, the following may be stated in terms
of the functions $(v_{i})_{i}$: find the coefficients
$\alpha_{i}^{n+1}$ in $u^{n+1}= \sum_{i=-1}^{N}\alpha_{i}^{n+1}
v_{i}$ such that
\begin{align}\label{eq14}
&\sum_{i=-1}^{N} \alpha_{i}^{n+1} \int_{\Omega} v_{i}v_{j}\, dx
  + \tau \sum_{i=-1}^{N} \alpha_{i}^{n+1} \int_{\Omega}k(u^{n}) \nabla v_{i}
  \nabla v_{j}\, dx +\tau \int_{\partial \Omega} \beta u^{n+1}  v_{j} \, ds \nonumber \\
 & = \sum_{i=-1}^{N} \alpha_{i}^{n} \int_{\Omega} v_{i}v_{j}\, dx
  + \tau \int_{\Omega} \sigma (u^{n}) |\nabla
\varphi^{n}|^{2} v_{j} \, dx.
\end{align}
In matrix notation, this may be expressed as
$$
\left (A+ \tau B \right) \alpha^{n+1}= f^{n}= f(n\tau), $$
where $$A= (a_{ij}) \mbox{ with element }a_{ij}= \int_{\Omega}
v_{i}v_{j}\, dx \, ,
$$
$$
B= (b_{ij}) \mbox{ with  } b_{ij}= \int_{\Omega}
k(u^{n}) \nabla v_{i} \nabla v_{j}\, dx \, ,$$ and
$$ \alpha^{n+1} \mbox { is the vector of unknows }
(\alpha_{i}^{n+1})_{i=-1}^{N}.
$$
Since the matrix $A$ and $B$ are Gram matrices, in particular they
are positive definite and invertible. Thus, the above system of
ordinary differential equations has obviously a unique solution.
We solve the system \eqref{eq11} for each time level. Estimating
each term of \eqref{eq11} separately, we have:
\begin{equation*}
\begin{split}
&\sum_{i=-1}^{N} \alpha_{i}^{n+1} \int_{\Omega} v_{i}v_{j}\, dx=
\sum_{i=-1}^{N} \alpha_{i}^{n+1} \int_{0}^{1}
v_{i}v_{j}\, dx \\
&\quad = \alpha_{j-1}^{n+1} \int_{x_{j-1}}^{x_{j}} v_{j-1}v_{j}\, dx
 +\alpha_{j}^{n+1} \int_{x_{j-1}}^{x_{j+1}} v_{j}^{2}\, dx
 + \alpha_{j+1}^{n+1} \int_{x_{j}}^{x_{j+1}} v_{j}v_{j+1}\, dx\\
 & \qquad +\alpha_{j-1}^{n+1} \int_{x_{j-1}}^{x_{j}} v_{j-1}v_{j}\, dx
 +\alpha_{j}^{n+1}\left( \int_{x_{j-1}}^{x_{j}} v_{j}^{2}\, dx
+\int_{x_{j}}^{x_{j+1}} v_{j}^{2}\, dx  \right)
 + \alpha_{j+1}^{n+1} \int_{x_{j}}^{x_{j+1}} v_{j}v_{j+1}\, dx.
\end{split}
\end{equation*}
Using the expression of $v_{j-1}, v_{j}$ and $v_{j+1}$, we obtain
\begin{equation}\label{eq15}
\sum_{i=-1}^{N} \alpha_{i}^{n+1} \int_{\Omega} v_{i}v_{j}\, dx=
\frac{h}{6}\alpha_{j-1}^{n+1}+\frac{2h}{3}\alpha_{j}^{n+1}+
\frac{h}{6}\alpha_{j+1}^{n+1}.
\end{equation}
In the same way, we have
\begin{equation*}
\begin{split}
&\sum_{i=-1}^{N} \alpha_{i}^{n+1} \int_{\Omega}k(u^{n}) \nabla v_{i}
\nabla v_{j}\, dx= \sum_{i=-1}^{N} \alpha_{i}^{n+1}
\int_{\Omega}k(u^{n})   \frac{\partial v_{i}}{\partial x}
 \frac{\partial v_{j}}{\partial x} \, dx \\
 & =\alpha_{j-1}^{n+1} \int_{x_{j-1}}^{x_{j}}k(u^{n}) \frac{\partial v_{j-1}}
 {\partial x}
 \frac{\partial v_{j}}{\partial x} \, dx
+ \alpha_{j}^{n+1} \int_{x_{j-1}}^{x_{j+1}}k(u^{n})( \frac{\partial
v_{j}}{\partial x})^{2}
  \, dx \\
 & \qquad + \alpha_{j+1}^{n+1} \int_{x_{j}}^{x_{j+1}}k(u^{n}) \frac{\partial v_{j}}{\partial x}
 \frac{\partial v_{j+1}}{\partial x} \, dx,\\
 & =- \frac{\alpha_{j-1}^{n+1}}{h^{2}} \int_{x_{j-1}}^{x_{j}}
 k(u^{n}) \, dx
 + \frac{\alpha_{j}^{n+1}}{h^{2}} \int_{x_{j-1}}^{x_{j+1}}
 k(u^{n}) \, dx
 - \frac{\alpha_{j+1}^{n+1}}{h^{2}} \int_{x_{j}}^{x_{j+1}}
 k(u^{n}) \, dx \\
 & \simeq - \frac{\alpha_{j-1}^{n+1}}{2h}
 \left( k(u^{n}(x_{j}))+ k(u^{n}(x_{j-1})) \right)
 + \frac{\alpha_{j}^{n+1}}{h}
 \left( k(u^{n}(x_{j+1}))+ k(u^{n}(x_{j-1})) \right)\\
 &\qquad - \frac{\alpha_{j+1}^{n+1}}{2h}
 \left( k(u^{n}(x_{j+1}))+ k(u^{n}(x_{j})) \right)\\
 & \simeq - \frac{\alpha_{j-1}^{n+1}}{2h}
 \left( k(\alpha_{j}^{n})+ k(\alpha_{j-1}^{n}) \right)
 + \frac{\alpha_{j}^{n+1}}{h}
 \left( k(\alpha_{j+1}^{n})+ k(\alpha_{j-1}^{n}) \right)
 - \frac{\alpha_{j+1}^{n+1}}{2h}
 \left( k(\alpha_{j+1}^{n})+ k(\alpha_{j}^{n}) \right).
\end{split}
\end{equation*}
On other hand, we similarly have
\begin{equation*}
\int_{\Omega} u^{n}v_{j}= \sum_{i=-1}^{N}
\alpha_{i}^{n}\int_{\Omega} v_{i} v_{j} \, dx=
\frac{h}{6}\alpha_{j-1}^{n}+ \frac{2h}{3} \alpha_{j}^{n}+
\frac{h}{6}\alpha_{j+1}^{n},
\end{equation*}
\begin{equation*}
\begin{split}
\int_{\Omega}\sigma(u^{n}) |\varphi_{x}^{n}|^{2} v_{j} \, dx & =
\sum_{j=1}^{N-1} \int_{x_{j}}^{x_{j+1}} \sigma(u^{n})
|\varphi_{x}^{n}|^{2} v_{j}(x) \, dx \\
& \simeq \frac{h}{2}\sum_{j=1}^{N-1} \left (\sigma(u^{n}(x_{j+1}))
|\varphi_{x}^{n}(x_{j+1})|^{2} v_{j}(x_{j+1}) + \sigma(u^{n}(x_{j}))
|\varphi_{x}^{n}(x_{j})|^{2} v_{j}(x_{j}) \right )\\
& \simeq  \frac{\sigma(\alpha_{j}^{n})}{h} (-\mu_{j-1}^{n}+
\mu_{j}^{n}+\mu_{j+1}^{n})^{2}.
\end{split}
\end{equation*}
It also holds:
\begin{equation*}
\begin{split}
\beta \int_{\partial \Omega = \{ 0, 1\}}
u^{n+1} v_{j} &= \beta u^{n+1}(1) v_{j}(1)-  \beta u^{n+1}(0) v_{j}(0)\\
& = \beta \alpha_{N}^{n+1} v_{j}(1)-  \beta \alpha_{0}^{n+1} v_{j}(0) \\
\end{split}
\end{equation*}
\begin{equation*}
=  \begin{cases} -\beta \alpha_{0}^{n+1}
 & \mbox{ if } j= 0,\\
 0 &  \mbox{ if } j=1 \ldots N-2,\\
0
& \mbox{ if } j=N-1. \\
  \end{cases}
\end{equation*}
Using together \eqref{eq14} and \eqref{eq15}, we get a system of
$N-1$ linear algebraic equations
\begin{equation}
\label{eq16}
\begin{split}
& \left ( \frac{h}{6} -\frac{\tau}{2h}(k(\alpha_{j}^{n})
  +k(\alpha_{j-1}^{n}) ) \right )\alpha_{j-1}^{n+1}
  +\left ( \frac{2}{3}h
+\frac{\tau}{h}(k(\alpha_{j+1}^{n})+k(\alpha_{j-1}^{n}) ) \right)
\alpha_{j}^{n+1}\\
& \qquad + \left ( \frac{h}{6} -\frac{\tau}{2h}(k(\alpha_{j+1}^{n})
  +k(\alpha_{j}^{n}) )\right) \alpha_{j+1}^{n+1} -\tau \beta  \alpha_{0}^{n+1} v_{j}(0)\\
  & =
  \frac{h}{6}\alpha_{j-1}^{n} + \frac{2h}{3}\alpha_{j}^{n}+
  \frac{h}{6}\alpha_{j+1}^{n}
  + \frac{\tau}{h}\sigma(\alpha_{j}^{n}) (-\mu_{j-1}^{n}+
\mu_{j}^{n}+\mu_{j+1}^{n})^{2}.
   \end{split}
\end{equation}
Using the boundary conditions, we find
\begin{equation*}
\begin{split}
&\alpha_{-1}^{n+1}= \alpha_{1}^{n+1}+ \left(
\frac{h\beta}{k(\alpha_{0}^{n})}-1 \right) \alpha_{0}^{n+1},\\
&\alpha_{-1}^{n}= \alpha_{1}^{n}+ \left(
\frac{h\beta}{k(\alpha_{0}^{n-1})}-1 \right) \alpha_{0}^{n},\\
& \alpha_{N}^{n+1}= \frac{k(\alpha_{N}^{n})}{\beta h +
k(\alpha_{N}^{n})} \alpha_{N-1}^{n+1},\\
& \alpha_{N}^{n}= \frac{k(\alpha_{N}^{n-1})}{\beta h +
k(\alpha_{N}^{n-1})} \alpha_{N-1}^{n}.
 \end{split}
\end{equation*}
From the initial condition we get
$$
\alpha_{0}^{0}= \alpha_{N}^{0}= 0 \, .
$$
Let
$$ a=\left ( \frac{h}{6}
-\frac{\tau}{2h}(k(\alpha_{0}^{n})+k(\alpha_{-1}^{n}) ) \right),
$$
$$
b=\left ( \frac{2h}{3} +\frac{\tau}{h}(k(\alpha_{1}^{n})
  +k(\alpha_{-1}^{n}) )
\right),
$$
and
$$
c=\left ( \frac{h}{6} -\frac{\tau}{2h}(k(\alpha_{1}^{n})
  +k(\alpha_{0}^{n}) )
\right).
$$
Substituting in \eqref{eq16}, we obtain the following system of
equations:

for $j=0$,
\begin{multline}
\left(a\left(\frac{\beta h}{k(\alpha_{0}^{n})}-1\right)+b - \tau \beta
 \right) \alpha_{0}^{n+1} +(a+c)\alpha_{1}^{n+1}\\
= \frac{h}{2} \left(1+ \frac{h \beta}{3k(\alpha_{0}^{n-1})}\right)
 \alpha_{0}^{n} + \frac{h}{3}\alpha_{1}^{n}
+ \frac{\tau}{h} \sigma(\alpha_{0}^{n})(2\mu_{0}^{n}+h
\overline{\varphi}(0))^{2} \, ;
\end{multline}
for $j= 1, \ldots, N-2$,
\begin{equation*}
\begin{split}
& \left ( \frac{h}{6} -\frac{\tau}{2h}(k(\alpha_{j}^{n})
  +k(\alpha_{j-1}^{n}) ) \right )\alpha_{j-1}^{n+1}
  +\left ( \frac{2h}{3}
+\frac{\tau}{h}(k(\alpha_{j+1}^{n})+k(\alpha_{j-1}^{n}) ) \right)
\alpha_{j}^{n+1}\\
& \qquad + \left ( \frac{h}{6} -\frac{\tau}{2h}(k(\alpha_{j+1}^{n})
  +k(\alpha_{j}^{n}) )\right) \alpha_{j+1}^{n+1}  \\
 & = \frac{h}{6}\alpha_{j-1}^{n} + \frac{2h}{3}\alpha_{j}^{n}+
  \frac{h}{6}\alpha_{j+1}^{n}+ \frac{\tau}{h}\sigma(\alpha_{j}^{n})(-\mu_{j-1}^{n}+
\mu_{j}^{n}+\mu_{j+1}^{n})^{2} \, ;
\end{split}
\end{equation*}
for $j= N-1$,
\begin{multline*}
d \alpha_{N-2}^{n+1} +\left(e+\frac{k(\alpha_{N}^{n})}{\beta h
+k(\alpha_{N}^{n})}f\right) \alpha_{N-1}^{n+1} \\
= \frac{h}{6}\alpha_{N-2}^{n} + \frac{h}{6}\left(4+
\frac{k(\alpha_{N}^{n-1})}{\beta h +k(\alpha_{N}^{n-1})} \right)
\alpha_{N-1}^{n} +
\frac{\tau}{h}\sigma(\alpha_{N-1}^{n})(2\mu_{N-1}^{n}-
\mu_{N-2}^{n}+h \overline{\varphi}(1))^{2},
\end{multline*}
where
$$
d= \left ( \frac{h}{6} -\frac{\tau}{2h}(k(\alpha_{N-1}^{n})
+k(\alpha_{N-2}^{n}) ) \right ),
$$
$$
e= \left ( \frac{2h}{3}
+\frac{\tau}{h}(k(\alpha_{N}^{n})+k(\alpha_{N-2}^{n}) \right),
$$
$$
f=\left ( \frac{h}{6} -\frac{\tau}{2h}(k(\alpha_{N}^{n})
  +k(\alpha_{N-1}^{n}) ) \right ).
$$


\section{An example}
\label{sec4}

In this section we give an example of a model of the thermistor
problem:
\begin{equation}
\label{eq:dif}
\left\{
\begin{gathered}
u_{t} = u_{xx} + \gamma  |\varphi_{x}|^{2} \\
\frac{\partial u}{\partial x} = - \beta u \mbox{ on } \partial
\Omega  \\
u(x, 0)= 0 , \quad 0 < x < 1
\end{gathered}
\right.
\end{equation}
\begin{equation}
\label{eq:elect}
\left\{
\begin{gathered}
(\sigma(u)\varphi_{x})_{x} = 0 \\
\frac{\partial \varphi}{\partial x} = 1 \mbox{ on } \partial
\Omega \\
\varphi (x, 0)= x, \quad 0\leq x  \leq 1.
\end{gathered}
\right.
\end{equation}

\begin{figure}
\begin{center}
\includegraphics[scale=0.65]{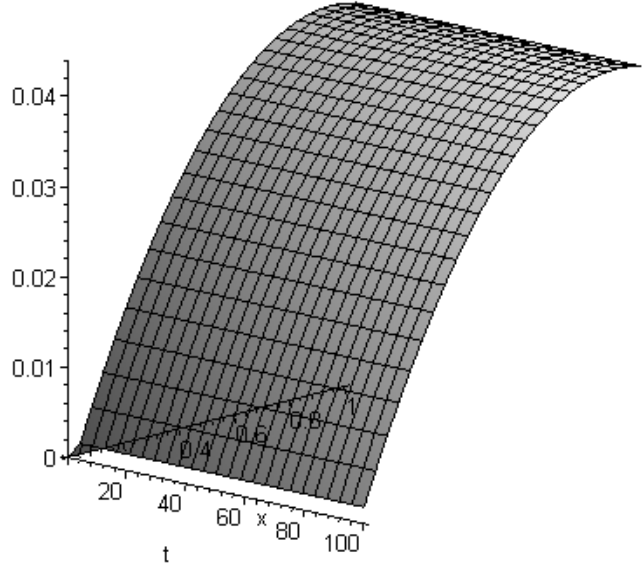}
\end{center}
\caption{The evolution of temperature.}
\label{fig:graph}
\end{figure}

The exact solution of the electrical potencial problem
\eqref{eq:elect} is $\varphi(t,x)=x$, $0 \le x \le 1$. Then, the
diffusion equation \eqref{eq:dif} can be reduced to the form
$$
\frac{\partial u}{\partial t}= \frac{\partial^{2} u}{\partial x^{2}}
+ \gamma.
$$
Using the proposed Galerkin finite element approach, we get the
following system of algebraic equations:

for $j=0$,
$$
\left(a_{1}(\beta h -1) +b_{1}- \tau \beta \right ) \alpha_{0}^{n+1}
+2a_{1}\alpha_{1}^{n+1}=\frac{h}{2}(1+ \frac{\beta
h}{3})\alpha_{0}^{n}+ \frac{h}{3} \alpha_{1}^{n}+ \gamma \tau h \, ;
$$
for $ j= 1, \ldots , N-2$,
\begin{equation*}
 a_{1} \alpha_{j-1}^{n+1}
  +b_{1}
\alpha_{j}^{n+1}
 + a_{1} \alpha_{j+1}^{n+1}
 = \frac{h}{6}\alpha_{j-1}^{n} + \frac{2h}{3}\alpha_{j}^{n}+
  \frac{h}{6}\alpha_{j+1}^{n}+ \gamma \tau h \, ;
\end{equation*}
for $j= N-1$,
\begin{equation*}
a_{1} \alpha_{N-2}^{n+1}
  +\left(b_{1}+\frac{a_{1}}{\beta h +1}\right) \alpha_{N-1}^{n+1}
 =
 \frac{h}{6}\alpha_{N-2}^{n} + \frac{h}{6}\left(4+
\frac{1}{1+\beta h } \right)
  \alpha_{N-1}^{n}
  +\gamma \tau h,
\end{equation*}
where
$$
a_{1}= \frac{h}{6}-\frac{\tau}{h} \, , \quad
b_{1} = \frac{2h}{3}+\frac{2\tau}{h}.
$$
We now show some results from numerical experiments performed using
our method and the computer algebra system \textsf{Maple~10}.
According with physical situations, we choose values of $\beta$ and
$\gamma$ verifying $\frac{1}{\beta}+ \frac{1}{2} \leq
\frac{1}{\gamma}$. In particular, we fixed $\beta=0.2$ and $\gamma =
0.1$. The calculation of the steady-state for the thermistor problem
is an important issue regarding the applications of the model in the
thermistor device. We obtained stable steady-state times for $\tau =
0.1, h=0.01$ (see Fig.~\ref{fig:graph}).

\section*{Acknowledgments}

The authors are grateful to the support of the \emph{Portuguese
Foundation for Science and Technology} (FCT) through the
\emph{Centre for Research in Optimization and Control} (CEOC) of
the University of Aveiro, cofinanced by the European Community
fund FEDER/POCI 2010, and the project SFRH/BPD/20934/2004.



\end{document}